\documentclass{article}
\usepackage[utf8]{inputenc}
\usepackage[margin=1 in]{geometry}
\usepackage{amsmath}
\usepackage{amsfonts}
\usepackage{amssymb}
\usepackage{amsthm}
\usepackage{mathrsfs}
\usepackage{enumitem}
\usepackage{lmodern}
\usepackage{hyperref}
\usepackage[vcentermath,enableskew]{youngtab}
\usepackage{ytableau}
\usepackage{graphbox}

\newtheorem{theorem}{Theorem}

\newtheorem{corollary}[theorem]{Corollary}
\newtheorem{conjecture}[theorem]{Conjecture}

\theoremstyle{definition}

\newtheorem{example}[theorem]{Example}

\newcommand{\mb}{\mathbb}
\newcommand{\tb}{\textbf}

\newcommand{\N}{\mathbb{N}}
\newcommand{\Z}{\mathbb{Z}}

\newcommand{\B}{\mathcal{B}}

\newcommand{\mc}{\mathcal}

\newcommand{\tn}{\textnormal}
\newcommand{\se}{\subseteq}

\newcommand{\lam}{\lambda}

\renewcommand{\B}{\boldsymbol}

\newcommand{\on}{\operatorname}

\newcommand{\mf}{\mathfrak}

\newcommand{\tcb}{\textcolor{blue}}

\newcommand{\im}{\includegraphics[width=1cm,align=c]} 

\usepackage[backend=bibtex,maxbibnames=10,citestyle=alphabetic]{biblatex}
\title{On the chromatic symmetric homology for star graphs}
\author{Laura Pierson \\ University of Waterloo \\ \href{mailto:lcpierson73@gmail.com}{lcpierson73@gmail.com}}

\begin{document}

\maketitle

\begin{abstract}
    The \emph{\tb{\tcb{chromatic symmetric function}}} $X_G$ is a power series that encodes the proper colorings of a graph $G$ by assigning a variable to each color and a monomial to each coloring such that the power of a variable in a monomial is the number of times the corresponding color is used in the corresponding coloring. The \emph{\tb{\tcb{chromatic symmetric homology}}} $H_*(G)$ is a doubly graded family of $\mb{C}[\mf{S}_n]$-modules that was defined by Sazdanovi\'c and Yip (2018) as a \emph{\tb{\tcb{categorification}}} of $X_G$. Chandler, Sazdanovi\'c, Stella, and Yip (2023) proved that $H_*(G)$ is a strictly stronger graph invariant than $X_G$, and they also computed or conjectured formulas for it in a number of special cases. We prove and extend some of their conjectured formulas for the case of \emph{\tb{\tcb{star graphs}}}, where one central vertex is connected to all other vertices and no other pairs of vertices are connected.
\end{abstract}

\section{Introduction}

The \emph{\tb{\tcb{chromatic symmetric function}}} $X_G$ was introduced by Stanley in \cite{stanley1995symmetric} as a symmetric function generalization of the chromatic polynomial \cite{birkhoff1913reducibility}. The \emph{\tb{\tcb{chromatic symmetric homology}}} $H_*(G)$ was introduced by Sazadnovi\'c and Yip in \cite{sazdanovic2018categorification} as a \emph{\tb{\tcb{categorification}}} of the chromatic symmetric function. The idea of categorification is to take a ``set-like" object and replace it with a ``category-like" object that has additional structure. Specifically, this often involves taking a polynomial and replacing it with a collection of modules over some ring, such as replacing \emph{\tb{\tcb{Schur polynomials}}} $s_\lam$ with \emph{\tb{\tcb{Specht modules}}} $\mc{S}_\lam$, which are irreducible representations of the symmetric group $\mf{S}_n,$ or equivalently, modules over the ring $\mb{C}[\mf{S}_n]$. The way these categorifications are constructed is generally analogous to the construction of the homology groups of a topological space. 

The idea of categorifying polynomials in this way began with the \emph{\tb{\tcb{Khovanov homology}}}, which was introduced by Khovanov in \cite{khovanov2000categorification} as a categorification of the \emph{\tb{\tcb{Jones polynomial}}} from knot theory. Sazdanovi\'c and Yip constructed the chromatic symmetric homology in \cite{sazdanovic2018categorification} as a categorification of the chromatic symmetric function. A corresponding vertex-weighted analogue of the chromatic symmetric homology was introduced in \cite{ciliberti2024deletion}. For a ring $R$, the chromatic symmetric homology $H_*(G;R)$ is constructed as $$H_*(G;R) = \bigoplus_{0\le i\le j \le |E(G)|} H_{i,j}(G;R),$$ where each $H_{i,j}(G;R)$ is an $R[\mf{S}_n]$-module. The chromatic symmetric homology $H_*(G)$ is a strictly stronger graph invariant than $X_G$ in the sense that there exist pairs of graphs with $X_{G_1} = X_{G_2}$ but $H_*(G_1) \ne H_*(G_2)$, as proven by Chandler, Sazdanovi\'c, Stella, and Yip in \cite{chandler2023strength}. In particular, the authors of \cite{chandler2023strength} conjectured that $H_*(G;\mb{Z})$ can detect planarity based on whether $H_{1,0}(G;\mb{Z})$ has $\Z_2$-torsion (i.e. contains a subgroup isomorphic to $\Z_2$). One direction of this conjecture was proven in \cite{ciliberti2021chromatic}, namely, that all nonplanar graphs contain $\Z_2$-torsion, but the other direction remains open.

The authors of \cite{chandler2023strength} also present a number of computations and conjectures about special cases of $H_*(G)$. In this chapter, we focus on the case of \emph{\tb{\tcb{star graphs}}} (graphs where one central vertex is connected to all others and no other pairs of vertices are connected). The authors of \cite{chandler2023strength} conjectured the following formulas involving the chromatic symmetric homology for star graphs:

\begin{conjecture}[\cite{chandler2023strength}]\label{conj:(n-2,2)}
    For $G$ the $n$-vertex star, the Specht module $\mc{S}_{(n-2, 2)}$ has multiplicity 1 in $H_{1,0}(G)$.
\end{conjecture}

\begin{conjecture}[\cite{chandler2023strength}]\label{conj:2^21^(n-4)}
    For $G$ the $n$-vertex star, the multiplicity of $\mc{S}_{2^21^{n-4}}$ in $H_{1,0}(G)$ is $\dbinom{n-2}2.$
\end{conjecture}

They also gave the following table (\cite{chandler2023strength}, \S 4.3) with their conjectures for the full values of $H_{1,0}(G)$ where $G$ is a star on 4, 5, 6, or 7 vertices:
\begin{table}[h!]
    \centering
    \renewcommand{\arraystretch}{1.5}
    \begin{tabular}{|c|c|}
        \hline
        $\B{G}$ & $\B{H_{1,0}(G;\mb{C})}$ \\ \hline
        \im{figs/star4} & $\mc{S}_{2^2}$ \\ \hline
        \im{figs/star5} & $\mc{S}_{2^21}^{\oplus3} \oplus \mc{S}_{32}$ \\ \hline
        \im{figs/star6} & $\mc{S}_{2^21^2}^{\oplus6} \oplus \mc{S}_{2^3}^{\oplus5} \oplus \mc{S}_{321}^{\oplus4} \oplus \mc{S}_{42}$ \\ \hline
        \im{figs/star7} & $\mc{S}_{2^21^3}^{\oplus10} \oplus \mc{S}_{2^31}^{\oplus16}\oplus \mc{S}_{321^2}^{\oplus 10} \oplus \mc{S}_{32^2}^{\oplus 9} \oplus \mc{S}_{421}^{\oplus5} \oplus \mc{S}_{52}$ \\ \hline
    \end{tabular}
    \label{tab:homology}
    \caption{The degree 0 chromatic symmetric homology for stars on 4 through 7 vertices.}
\end{table}

\bigskip

We verify Conjectures \ref{conj:(n-2,2)} and \ref{conj:2^21^(n-4)} as well as the conjectured values from Table \ref{tab:homology} as special cases of the following result:

\begin{theorem}\label{thm:general_l_k}
    For $G$ the $n$-vertex star, the multiplicity of $\mc{S}_{\ell 2^k 1^{n-\ell - 2k}}$ in $H_{1,0}(G)$ is $$\binom{n-1}{\ell-1}f^{2^k1^{n-\ell-2k}} - f^{\ell 2^k 1^{n-\ell-2k}},$$ where $f^\lam$ is the number of standard Young tableaux of shape $\lam.$
\end{theorem}

In particular, Conjectures \ref{conj:(n-2,2)} and \ref{conj:2^21^(n-4)} follow immediately from the following special case of Theorem \ref{thm:general_l_k}:

\begin{corollary}\label{cor:l21^(n-l-2)}
    For $G$ the $n$-vertex star multiplicity of $\mc{S}_{\ell 2 1^{n-\ell-2}}$ in $H_{1,0}(G)$ is $\dbinom{n-2}\ell$.
\end{corollary}

We also give a more explicit formula for the special case where all parts of $\lam$ have size 1 or 2, which we will use to verify some of the values in Table \ref{tab:homology}:

\begin{corollary}\label{cor:2^k1^(n-2k)}
    For $G$ the $n$-vertex star, the multiplicity of $\mc{S}_{2^k1^{n-2k}}$ in $H_{1,0}(G)$ is $$(n-2k+1)\left(\binom{n-1}{k-1}-\frac1k\binom n{k-1}\right).$$
\end{corollary}


Finally, we conjecture that Theorem \ref{thm:general_l_k} gives essentially a full description of the degree 0 chromatic symmetric homology for star graphs, in the sense that all other multiplicities are 0 (except $H_{0,0}(G) = \mc{S}_{1^n}$):

\begin{conjecture}
    Whenever $i\ge 2$ or $\lam_2 \ge 3,$ the multiplicity of $\mc{S}_\lam$ in $H_{i,0}(G)$ is 0.
\end{conjecture}
    
The remainder of this paper is organized as follows. We first introduce some relevant background on graphs and symmetric functions in \S \ref{sec:background}. We then give the general construction of $H_*(G)$ in \S \ref{sec:homology_background}. In \S\ref{sec:deg_0} and \S \ref{sec:decomposition}, we explain how to concretely compute the degree 
$j=0$ case $H_{i,0}(G;\mb{C})$ and its decomposition into Specht modules using the combinatorics of Young tableaux. In \S \ref{sec:general_l_k} we prove Theorem \ref{thm:general_l_k}. In \S \ref{sec:l21^(n-l-2)}, we prove Corollary \ref{cor:l21^(n-l-2)}, and in \S \ref{sec:2^k1^(n-2k)}, we prove Corollary \ref{cor:2^k1^(n-2k)}. Finally, in \S \ref{sec:table_vals}, we use our results to verify the conjectured values from Table \ref{tab:homology}.

\section{Background on graphs and symmetric functions}\label{sec:background}

A \emph{\tb{\tcb{graph}}} $G$ consists of a set $V(G)$ of \emph{\tb{\tcb{vertices}}} and a set $E(G)$ of unordered pairs $uv$ of vertices, called \emph{\tb{\tcb{edges}}}. The \emph{\tb{\tcb{star graph}}} on vertex set $\{1,2,\dots,n\}$ is the graph with edges $1i$ for each $i\ne 1$ and no other edges. Write $\N = \{1,2,3,\dots\}$ for the set of natural numbers. A \emph{\tb{\tcb{proper coloring}}} of $G$ is a function $\kappa:V(G)\to \N$ such that $\kappa(u)\ne \kappa(v)$ whenever $uv \in E(G),$ which we can think of as an assignment of a color to each vertex such that no two adjacent vertices get the same color.  A \emph{\tb{\tcb{symmetric function}}} is a polynomial $f(x_1,x_2,\dots)$ that stays the same under any permutation of the variables, i.e. $f(x_1,x_2,\dots) = f(x_{\sigma(1)},x_{\sigma(2)},\dots)$ for any permutation $\sigma$ of $\N.$ The \emph{\tb{\tcb{chromatic symmetric function}}} is $$X_G(x_1,x_2\dots) := \sum_\kappa \prod_{v\in V(G)} x_{\kappa(v)}.$$

A \emph{\tb{\tcb{partition}}} $\lam = (\lam_1,\dots,\lam_\ell) = \lam_1\dots \lam_\ell$ is a finite nondecreasing sequence of positive integers, and the $\lam_i$'s are its \emph{\tb{\tcb{parts}}}. We also write $\lam = i_1^{j_1}i_2^{j_2}\dots i_\ell^{j_\ell}$ to mean the partition that has $j_k$ parts of size $i_k$ for each $k=1,2,\dots,\ell.$ The \emph{\tb{\tcb{Young diagram}}} for $\lam$ consists of $\ell$ left-justified rows of boxes such that the $i$th row from the top has $\lam_i$ boxes. A \emph{\tb{\tcb{Young tableau}}} is a filling of the boxes of the Young diagram with positive integers. Its \emph{\tb{\tcb{content}}} is the sequence $(\mu_1,\mu_2,\dots)$ such that $\mu_i$ is the number of entries equal to $i.$ A tableau is a \emph{\tb{\tcb{semistandard Young tableau (SSYT)}}} if the entries increase weakly across rows and strictly down columns. The \emph{\tb{\tcb{Kostka number}}} $K_{\lam\mu}$ counts the number of SSYT of shape $\lam$ and content $\mu$. An SSYT is a \emph{\tb{\tcb{standard Young tableau (SYT)}}} if the content is $1^n$, i.e. the entries are some permutation of $\{1,2,\dots,n\}.$ We write $$f^{\lam} := K_{\lam,1^n} = |{\sf{SYT}}(\lam)|$$ for the number of SYT of shape $\lam.$

\begin{example}\label{ex:SYT}
    We have $f^{32} = 5,$ with the 5 elements of ${\sf{SYT}}(32)$ shown below:
    \begin{center}
        \young(123,45) \hspace{1cm} \young(124,35) \hspace{1cm} \young(125,34) \hspace{1cm} \young(134,25) \hspace{1cm} \young(135,24).
    \end{center}
\end{example}

The number of SYT of shape $\lam$ can be calculated using the \emph{\tb{\tcb{hook length formula}}}. For each box $b \in \lam$ in the Young diagram of shape $\lam$, the \emph{\tb{\tcb{hook}}} ${\sf{hook}}(b)$ consists of $b$ itself together with all the boxes below it in its column or to the right of it in its row (forming a $\Gamma$ shape). Then the formula is $$f^\lam = \frac{|\lam|!}{\prod_{b\in \lam} |{\sf{hook}}(b)|},$$ where $|{\sf{hook}}(b)|$ denotes the number of boxes in ${\sf{hook}}(b).$

\begin{example}
    For $\lam = 32,$ the hooks are shown below:
    $$ \ytableausetup{boxsize=1.2em}
    \ydiagram[*(green)]{3,1}*[*(white)]{3,2} \hspace{1cm}
    \ydiagram[*(green)]{1+2,1+1}*[*(white)]{3,2} \hspace{1cm}
    \ydiagram[*(green)]{2+1}*[*(white)]{3,2} \hspace{1cm}
    \ydiagram[*(green)]{3+0,2}*[*(white)]{3,2} \hspace{1cm}
    \ydiagram[*(green)]{3+0,1+1}*[*(white)]{3,2}.$$
    Their lengths are 4, 3, 1, 2, and 1, respectively, so the hook length formula gives $$f^{32} = \frac{5!}{4\cdot 3 \cdot 1 \cdot 2 \cdot 1} = 5,$$ which matches what we found in Example \ref{ex:SYT}.
\end{example}

The \emph{\tb{\tcb{symmetric group}}} $\mf{S}_n$ is the group of permutations of the set $\{1,2,\dots,n\}.$ A \emph{\tb{\tcb{representation}}} of a group $G$ is a linear action of $G$ on some vector space, or equivalently, a module $V$ over the ring $\mb{C}[G]$. A representation $V$ is \emph{\tb{\tcb{irreducible}}} if it has no proper submodule $W \se V$ such that $\sigma W \se W$ for all $\sigma \in G.$ The irreducible representations of $\mf{S}_n$ are called \emph{\tb{\tcb{Specht modules}}}, and are indexed by the partitions $\lam$ of $n.$ Given a representation $V$ of a subgroup $H \se G,$ the corresponding \emph{\tb{\tcb{induced representation}}} of $G$ is $\on{Ind}_H^G(V) := \mb{C}[G] \otimes_{\mb{C}[H]} V.$

\section{General construction of $H_*(G)$}\label{sec:homology_background}

Let $n := |V(G)|$ and $m := |E(G)|.$ For a ring $R$ (generally $\mb{C}$ or $\mb{Z}$), the chromatic symmetric homology $H_*(G;R)$ is defined by building a sequence of \emph{\tb{\tcb{chain modules}}} $C_i(G)$ with associated \emph{\tb{\tcb{differentials}}} or \emph{\tb{\tcb{chain maps}}} $d_i:C_i(G)\to C_{i-1}(G)$, where $0 \le i \le m$: $$0 \overset{d_{m+1}}{\longrightarrow} C_m(G) \overset{d_m}{\longrightarrow} C_{m-1}(G) \overset{d_{m-1}}{\longrightarrow} \dots \overset{d_2}{\longrightarrow} C_1(G) \overset{d_1}{\longrightarrow} C_0(G) \overset{d_0}{\longrightarrow} 0.$$ Each $C_i(G)$ is a representation of $\mf{S}_n,$ or equivalently, an $R[\mf{S}_n]$-module, and each $d_i$ is a homomorphism of $R[\mf{S}_n]$-modules. The \emph{\tb{\tcb{homology groups}}} associated to $G$ are then defined by $$H_i(G;R) := \ker(d_i)/\on{im}(d_{i+1}).$$ Each chain module has an associated \emph{\tb{\tcb{grading}}} $C_i(G;R) = \bigoplus_{j=0}^i C_{i,j}(G;R),$ and the chain maps and homology groups inherit this grading, so we can write $d_i = \bigoplus_{j=0}^i d_{i,j}$ and $H_i(G;R) = \bigoplus_{j=0}^i H_{i,j}(G,R).$ The key connection to $X_G$ is that when $R = \mb{C}$, if we break each $H_{i,j}$ into a direct sum of Specht modules, replace each Specht module with the associated Schur polynomial, and sum the resulting Schur polynomials with the appropriate signs, we recover the chromatic symmetric function: $$H_{i,j}(G;\mb{C}) = \bigoplus_\lam \mc{S}_\lam^{\oplus n_{ij}^\lam} \hspace{0.5cm} \longrightarrow \hspace{0.5cm} X_G = \sum_{i,j}(-1)^{i+j}n_{ij}^\lam s_\lam,$$ where $n_{ij}^\lam$ is the multiplicity of $\mc{S}_\lam$ in $H_{i,j}(G;\mb{C}).$ Assume the vertices of $G$ are labeled 1 through $n.$ The chain module $C_i(G)$ can be constructed as a direct sum of $\mb{C}[\mf{S}_n]$-modules $\mc{M}_F$, one for each subset $F \se E$ with $|F| = i$: $$C_i(G) := \bigoplus_{F\se E,|F|=i}\mc{M}_F.$$ Each such subset $F$ corresponds to a \emph{\tb{\tcb{spanning subgraph}}} of $G$, where ``spanning" means we take the vertex set of the subgraph to be the full vertex set $V(G).$ If $b_1,\dots,b_r$ are the sizes of the connected components of $F$, then the definition of $\mc{M}_F$ is $$\mc{M}_F := \on{Ind}_{\mf{S}_{b_1}\times \dots \times \mf{S}_{b_r}}^{\mf{S}_n}(\mc{L}_{b_1} \otimes \dots \otimes \mc{L}_{b_r}),$$ where $\on{Ind}_H^G(V)$ denotes the \emph{\tb{\tcb{induced representation}}} of $G$ corresponding to a representation $V$ of $H,$ $$\mc{L}_b := \bigoplus_{j=0}^{b-1} \mc{S}_{(b-j,1^j)},$$ and $\mc{S}_\lam$ is the \emph{\tb{\tcb{Specht module}}} corresponding to the partition $\lam.$ The \emph{\tb{\tcb{degree (or grading)}}} is defined such that each $\mc{S}_{(b-j,1^j})$ piece has degree $j,$ and the degree of a tensor product of two representations is the sum of their degrees. In particular, the degree 0 parts are $$(\mc{M}_F)_0 = \on{Ind}_{\mf{S}_{b_1} \times \dots \times \mf{S}_{b_r}}^{\mf{S}_n}(\mc{S}_{b_1} \otimes \dots \otimes \mc{S}_{b_r}).$$ 

Then for each pair $F$ and $F'$ such that $F'$ is formed by removing one edge from $F$, the differential map $d_{F,F'}$ is the natural embedding of $\mc{M}_F$ into $\mc{M}_{F'}.$ (The precise construction of this unique ``natural embedding" is a bit technical and is given in \cite{sazdanovic2018categorification}, \S 2.3.)

Finally, given an ordering on $|E(G)|$, the differential maps $d_i$ are defined by $$d_i := \bigoplus_{|F|=i,|F'|=i-1} (-1)^{\tn{\# of edges in $F$ that come after the removed edge}}\cdot d_{F,F'}.$$ It can be shown that these maps are in fact differentials in the sense that $d_{i}\circ d_{i+1} = 0,$ or equivalently, $\on{im}(d_{i+1}) \se \ker(d_i).$

\section{The degree $j=0$ case $H_{i,0}(G)$}\label{sec:deg_0}

Our focus from now on will be on the case $j=0$ and $R=\mb{C},$ so will now describe how $H_{i,0}(G)$ can be explicitly computed using the combinatorics of Young tableaux.

In that case, the $\mb{C}[\mc{S}_n]$-modules $\mc{M}_F$ (which we henceforth use as shorthand for $(\mc{M}_F)_0$) can alternatively be constructed as $$\mc{M}_F := \mb{C}[\mf{S}_n]\cdot a_{T(F)},$$ where $\mb{C}[\mf{S}_n]\cdot a_{T(F)} := \{\sigma \cdot a_{T(F)} : \sigma \in \mb{C}[\mf{S}_n]\},$ and $a_{T(F)}$ is a particular element of $\mb{C}[\mf{S}_n]$ that we will define momentarily. To define $a_{T(F)},$ let $\lam(F)$ be \emph{\tb{\tcb{partition type}}} of $F$, meaning the partition whose parts are the sizes of the connected components of $F,$ arranged in decreasing order. We construct a Young tableau $T(F)$ of shape $\lam(F)$ as follows. For each connected component of $F$, list the vertices of $F$ in increasing order. Order these lists first in decreasing order by length, and then within that in increasing order by the smallest element. Then use the $i$th list to fill the $i$th row of the Young diagram. Note that the resulting Young tableau need not be standard, as the numbers increase across each row but may not increase down each column.

\begin{example}\label{ex:graph_tableau}
    Let $F$ be the spanning subgraph shown below:
    \begin{center}
        \includegraphics[width=3cm]{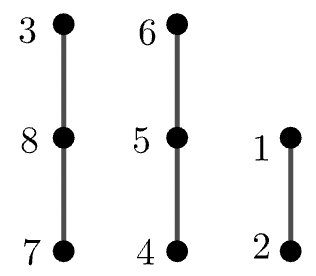}
    \end{center}
    Then $\lam(F) = 332,$ and $$T(F) = \young(378,456,12).$$
\end{example}

The \emph{\tb{\tcb{row symmetrizer}}} $a_T \in \mb{C}[\mf{S}_n]$ for a tableau $T$ is the sum of all $\sigma \in\mf{S}_n$ such that $\sigma$ permutes the elements in each row of $T(F)$ among themselves, and the associated \emph{\tb{\tcb{permutation module}}} is $$\mc{M}_T := \mb{C}[\mf{S}_n]\cdot a_T.$$ We can factor $a_T$ as the product over rows of the sum of all permutations of elements in that row.

\begin{example}
    Write $e$ for the identity permutation and $(i_1i_2 \dots i_k)$ for the $k$ cycle that maps $i_j$ to $i_{j+1}$ for $1\le j\le k-1$ and $i_k$ to $i_1.$ Then for the $F$ from Example \ref{ex:graph_tableau}, $$a_{T(F)} = (e + (37) + (38) + (78) + (378) + (387))\cdot (e + (45) + (46) + (56) + (456) + (465)) \cdot (e + (12)).$$
\end{example}

\section{Decomposition of $H_{i,0}(G)$ into Specht modules}\label{sec:decomposition}

The \emph{\tb{\tcb{column symmetrizer}}} $b_T\in \mb{C}[\mf{S}_n]$ is the sum of $\on{sign}(\sigma)\cdot \sigma$ over all permutations $\sigma$ that permute the elements of each column among themselves.

\begin{example}
    For the $F$ from Example \ref{ex:graph_tableau}, $$b_{T(F)} = (e - (13) - (14) - (34) + (134) + (143))\cdot (e - (25) - (27) - (57) + (257) + (275)) \cdot (e - (68)).$$
\end{example}

The \emph{\tb{\tcb{Young symmetrizer}}} for $T$ is $c_T := b_T a_T.$ The corresponding \emph{\tb{\tcb{Specht module}}} is $$\mc{S}_T := \mb{C}[\mf{S}_n] \cdot c_T.$$ For any tableaux $S$ and $T$ of the same shape $\lam,$ $\mc{S}_T$ and $\mc{S}_{T'}$ are isomorphic, so we can write $\mc{S}_\lam$ for this isomorphism class of Specht modules. Similarly, $\mc{M}_T$ and $\mc{M}_{T'}$ are isomorphic when $T$ and $T'$ have the same shape $\mu,$ so we can write $\mc{M}_\mu$ for this isomorphism class. 

It turns out that the Specht modules $\mc{S}_\lam$ are precisely the \emph{\tb{\tcb{irreducible representations}}} of $\mf{S}_n$ (see \S \ref{sec:rep_theory}), which means any $\mb{C}[\mf{S}_n]$ module can be written as a direct sum of Specht modules in a way that is unique up to isomorphism. In the case of $\mc{M}_F,$ this decomposition turns out to be $$\mc{M}_F \cong \bigoplus_\lam \mc{S}_\lam^{\oplus K_{\lam,\lam(F)}},$$ where $K_{\lam\mu}$ is the \emph{\tb{\tcb{Kostka number}}} that counts the number of SSYT of shape $\lam$ and content $\mu.$

We can realize this decomposition for $\mc{M}_F \cong \mc{M}_{\lam(F)}$ as $$\mc{M}_F = \bigoplus_\lam \bigoplus_{j=1}^{K_{\lam,\lam(F)}}\mb{C}[\mf{S}_n]\cdot v_{X_F^j}^{T(\lam)}.$$ It remains to explain what $x_F^j$, $T(\lam),$ and $v_{X_F^j}^{T(\lam)}$ are. For $S,T \in {\sf{SYT}}(\lam),$ write $\sigma_{S,T} \in \mf{S}_n$ for the permutation such that $\sigma_{S,T}\cdot S = T.$ Then $$v_S^T := c_T \cdot \sigma_{S,T} = \sigma_{S,T} \cdot c_S.$$ The tableau $T(\lam) \in {\sf{SYT}}(\lam)$ is constructed by filling in the numbers $1,2,\dots,n$ in order row by row. For instance, for $\lam = 432,$ $$T(\lam) = \young(1234,567,89).$$ Finally, the tableaux $X_F^j$ are constructed as follows. Let $T_1,\dots,T_{K_{\lam,\lam(F)}}$ be the SSYT of shape $\lam$ and content $\lam(F).$ Then we turn $T_j$ into the corresponding tableau $X_F^j$ of shape $\lam$ by replacing the 1's in $T_j$ with the entries from the first row of $T(F)$ in order, the 2's with the entries from the second row of $T(F)$ in order, and so on, and then reordering so the entries in each row are increasing. The resulting tableau $X_F^j$ will have entries 1 through $n$ each showing up exactly once, and the entries will increase across each row, but it will not necessarily be standard as the entries need not increase going down the columns.

\begin{example}
    Let $F$ be the spanning subgraph from Example \ref{ex:graph_tableau}, so $\lam(F) = 332,$ and let $\lam = 431.$ Then $K_{\lam,\lam(F)} = 2,$ so there are two copies of $\mc{S}_\lam$ in $\mc{M}_T.$ The two the SSYT of shape $\lam$ and content $\lam(F)$ are $$T_1 = \young(1112,223,3), \hspace{1cm} T_2 = \young(1113,222,3).$$ To get the tableaux $X_F^1$ and $X_F^2,$ we replace the three 1's in each tableau by the 3, 7, and 8 from the first row of $\lam(F)$, we replace the three 2's by 4, 5, and 6 from the second row, and we replace the two 3's by 1 and 2 from the third row. This gives: $$\young(3784,561,2), \hspace{1cm} \young(3781,456,2).$$ Then we reorder each row so the entries are in increasing order to get $$X_F^1 = \young(3478,156,2), \hspace{1cm} X_F^2 = \young(1378,456,2).$$ Thus, the two copies of $\mc{S}_\lam$ in $\mc{M}_T$ are given by $$\mc{M}_T|_{\mc{S}_{431}} = \mb{C}[\mf{S}_n]\cdot v_{\young(3478,156,2)}^{T(\lam)} \oplus \mb{C}[\mf{S}_n]\cdot v_{\young(3781,456,2)}^{T(\lam)},$$ where $$T(\lam) = \young(1378,456,2).$$ Note that neither $X_F^1$ not $X_F^2$ is an SYT, since the entries do not always increase down the columns.
\end{example}

\section{Proof of Theorem \ref{thm:general_l_k}}\label{sec:general_l_k}

Let $\lam = \ell 2^k 1^{n-\ell-2k}.$ The restriction of $C_0(G)$ to $\mc{S}_\lam$ is $$C_0(G)|_{\mc{S}_\lam} = \bigoplus_{Y_i\in {\sf{SYT}}(\lam)} \mb{C}[\mf{S}_n]\cdot v_{Y_i}^{T(\lam)} \se \mb{C}[\mf{S}_n],$$ so the multiplicity of $\mc{S}_\lam$ in $C_0(G)$ is $|{\sf{SYT}}(\lam)| = f^\lam.$

To find $C_1(G)|_{\mc{S}_\lam},$ we need to consider the spanning subgraphs of $G$ with a single edge. Label the vertices of $G$ as $1,\dots,n,$ and without loss of generality assume vertex 1 is the central vertex. Then there are $n-1$ spanning subgraphs of $G$ with one edge, each of the form $F_{1i}$ for some $i\ne 1.$ Since $\lam(F_{1i}) = 21^{n-2}$ for all $i$, the restriction of $C_1(G)$ to $\mc{S}_\lam$ is $$C_1(G)|_{\mc{S}_\lam} = \bigoplus_{i=2}^n \mc{M}_{F_{1i}}|_{\mc{S}_\lam} = \bigoplus_{i=2}^n \bigoplus_{j=1}^{K_{\lam,21^{n-2}}} \mb{C}[\mf{S}_n]\cdot v_{X_i^j}^{T(\lam)},$$ where each tableau $X_i^j$ is formed by taking an SSYT of shape $\lam$ and content $\lam(F_{1i}) = 21^{n-2}$, then relabeling it so the second 1 (which is necessarily the second entry of the top row) becomes an $i$ and the entries greater than $i-1$ each get increased by 1, then reordering the first row so the entries are in increasing order. Every $Y_k$ is equal to $X_i^j$ for some choice of $i$ and $j$ (e.g. by taking $i$ to be the second entry in the top row of $Y_k$), so $d_1|_{\mc{S}_\lam}$ is surjective. Thus, by rank-nullity, the multiplicity of $\mc{S}_\lam$ in $\ker(d_1)$ is equal to its multiplicity in $C_1(G)$ minus its multiplicity in $C_0(G).$ The multiplicity in $C_0(G)$ is simply $f^\lam.$ We claim that the multiplicity in $C_1(G)$ is $$(\ell-1)\cdot\binom{n-1}{\ell-1}f^{2^k 1^{n-\ell-k}}.$$ To see this, note that to form an $X_i^j,$ we can choose any $\ell-1$ numbers besides 1 to fill the first row (in $\binom{n-1}{\ell-1}$ ways), choose any one of them to be $i$ (in $\ell-1$ ways), and then fill the remaining rows by taking any SSYT of shape $2^k 1^{n-\ell-k}$ (in $f^{2^k1^{n-\ell-2k}}$ ways) and relabeling it so the 1 becomes the smallest remaining entry, the 2 becomes the second smallest, and so on.

For $C_2(G)|_{\mc{S}_\lam},$ we consider the $\binom{n-1}2$ spanning subgraphs of $G$ with two edges, which are each of the form $F_{1i,1j}$ for $i$ and $j.$ Every one of these spanning subgraphs has partition type $\lam(F_{1i,1j}) = 31^{n-3},$ so $$C_1(G)|_{\mc{S}_\lam} = \bigoplus_{2\le i<j \le n} \mc{M}_{F_{1i,1j}}|_{\mc{S}_\lam} = \bigoplus_{2\le i<j \le n} \bigoplus_{k=1}^{K_{\lam,31^{n-3}}}\mb{C}[\mf{S}_n]\cdot v_{W_{ij}^k}^{T(\lam)},$$ where each $W_{ij}^k$ is formed from an SSYT of shape $\lam$ and content $31^{n-3}$ by replacing two of the 1's with $i$ and $j,$ adding 1 to each entry between $i$ and $j-2$, adding 2 to each entry between $j-1$ and $n-2$, and reordering the first row so its entries are increasing. Every $W_{ij}^k$ formed in this manner is also of the form $X_i^r$ and an $X_j^s$ for some $r$ and $s$. Thus, the map $d_2|_{\mc{S}_\lam}$ directly embeds each $\mb{C}[\mf{S}_n]\cdot v_{W_{ij}^k}^{T(\lam)} \se \mc{M}_{1i,1j}$ into $$\mb{C}[\mf{S}_n]\cdot v_{X_i^r}^{T(\lam)} \oplus [\mf{S}_n]\cdot v_{X_j^s}^{T(\lam)} \se \mc{M}_{1i}\oplus \mc{M}_{1j},$$ via a map that looks like $$d_2|_{\mc{S}_\lam}(\tau_{ij}^k) = (-\tau_{ij}^k) \oplus \tau_{ij}^k$$ for each $\tau_{ij}^k \in \mb{C}[\mf{S}_n]\cdot v_{W_{ij}^k}^{T(\lam)}.$ We write the first component with a negative sign because can assume without loss of generality that the edge $1i$ comes before the edge $1j$ in the edge ordering, and hence the map $\mc{M}_{1i,1j} \to \mc{M}_{1i}$ has negative sign (since there is one edge in $F$, namely $1j$, that comes after $1i$ in the ordering) while the map $\mc{M}_{1i,1j} \to \mc{M}_{1j}$ has positive sign (since there are no edges in $F$ coming after $1j$ in the ordering).

Each tableau that shows up as $X_i^j$ for some $i$ actually shows up as $X_i^j$ for $\ell-1$ different values of $i$, since we can take $i$ to be any entry in the top row except the 1. It also shows up as $W_{ij}^k$ for $\binom{\ell-1}2$ pairs $(i,j)$, since we can take $i$ and $j$ two be any two entries in the top row except the two 1's. Furthermore, a tableau of shape $\lam$ is one of the $X_i^j$'s if and only if it is one of the $W_{ij}^k$'s, since the only requirement in both cases is that the first row be in increasing order and the portion below the first row be standard.

To find the multiplicity of $\mc{S}_\lam$ in $d_2,$ consider the restriction of $d_2$ to submodules of the form $\mb{C}[\mf{S}_n]\cdot v_S^{T(\lam)}$ for a particular tableau $S$. There are $\ell-1$ such modules in $C_1(G)$, one for each entry in the top row of $S$ besides the 1, each isomorphic to a copy of $\mc{S}_\lam.$ From the description of $d_2$ above, an element $\sigma_1 \oplus \dots \oplus \sigma_{\ell-1}$ in the direct sum of these $\ell-1$ submodules is in $\on{im}(d_2)$ if and only if $\sigma_1 + \dots + \sigma_{\ell-1} = 0.$ Thus, $\sigma_1,\dots,\sigma_{\ell-1}$ can be chosen independently, and $\sigma_{\ell-1}$ is then determined. Thus, for every $\ell-1$ copies of $\mc{S}_\lam$ in $C_1(G)$, there are $\ell-2$ copies in $\on{im}(d_2),$ so the multiplicity of $\mc{S}_\lam$ in $\on{im}(d_2)$ is $\frac{\ell-1}{\ell-2}$ times its multiplicity in $C_1(G),$ or $$(\ell-2)\cdot \binom{n-1}{\ell-1}f^{2^k1^{n-\ell-2k}}.$$ Then since $H_{1,0}(G) := \ker(d_1)/\on{im}(d_2),$ the multiplicity of $\mc{S}_\lam$ in $H_{1,0}(G)$ is equal to its multiplicity in $\ker(d_1)$ minus its multiplicity in $\on{im}(d_2),$ or $$(\ell-1)\cdot \binom{n-1}{\ell-1}f^{2^k1^{n-\ell-2k}} - f^\lam - (\ell-2)\cdot \binom{n-1}{\ell-1}f^{2^k1^{n-\ell-2k}} = \binom{n-1}{\ell-1}f^{2^k1^{n-\ell-2k}} - f^\lam,$$ as claimed. \qed

\section{Proof of Corollary \ref{cor:l21^(n-l-2)}}\label{sec:l21^(n-l-2)}

Let $\lam = \ell 2 1^{n-\ell-2}.$ By Theorem \ref{thm:general_l_k}, the multiplicity of $\mc{S}_{\ell 21^{n-\ell-2}}$ in $H_{1,0}(G)$ is $$\binom{n-1}{\ell-1}f^{21^{n-\ell-2}} - f^\lam.$$ We note that $$f^{21^{n-\ell-2}} = n-\ell-1,$$ since to get an SYT of shape $21^{n-\ell-2}$, we can choose any of the $n-\ell-1$ numbers from $2,\dots,n-\ell$ to be the second entry in the first row, and then the rest is determined since the entries in the first column must be in increasing order. To compute $f^\lam,$ we draw out the Young diagram of shape $\lam$ and label each box with its corresponding hook length:
$$\ytableausetup{boxsize=3em}
\begin{ytableau}
    n-1 & \ell & \ell-2 & \dots & 2 & 1 \\
    n-\ell & 1 \\
    \scriptstyle{n-\ell-2} \\
    \vdots \\
    2 \\
    1
\end{ytableau}$$
Thus, the hook length formula gives $$f^\lam = \frac{n!}{(n-1)(n-\ell)\ell\cdot(n-\ell-2)!\cdot(\ell-2)!} = \frac{n(\ell-1)}{n-\ell}\binom{n-2}\ell.$$ Putting this together, our multiplicity is $$(n-\ell-1)\binom{n-1}{\ell-1} - \frac{n(\ell-1)}{n-\ell}\binom{n-2}\ell = \frac{(n-1)\ell - n(\ell-1)}{n-\ell}\binom{n-2}\ell = \binom{n-2}\ell,$$ as claimed. Conjectures \ref{conj:(n-2,2)} and \ref{conj:2^21^(n-4)} then follow immediately as the special cases $\ell=n-2$ and $\ell=2.$ \qed

\section{Proof of Corollary \ref{cor:2^k1^(n-2k)}}\label{sec:2^k1^(n-2k)}

Let $\lam = 2^k 1^{n-2k}.$ Since $\ell=2,$ Theorem \ref{thm:general_l_k} implies that the multiplicity of $\mc{S}_\lam$ in $H_{1,0}(G)$ is $$(n-1)K_{\lam,21^{n-2}} - f^\lam.$$ By the hook-length formula, we get $$f^\lam = \frac{n!}{(n-k+1)\dots(n-2k+2)\cdot k! \cdot (n-2k)!} = \frac{n-2k+1}{k}\binom{n}{k-1},$$ since we can label each box with its hook length as shown below:
$$\begin{ytableau}
    \scriptstyle{n-k+1} & k \\
    n-k & k-1 \\
    \vdots & \vdots \\
    \scriptstyle{n-2k+2} & 1 \\
    n-2k \\
    \vdots \\
    1
\end{ytableau}$$
For $K_{\lam,21^{n-2}}$, note that to get an SSYT of shape $\lam$ and content $21^{n-2},$ the two 1's must be in the top row, and then the number of ways to fill the remaining part is equal to the number of SYT of shape $2^{k-2}1^{n-2k}.$ Using the formula above for $f^{2^k1^{n-2k}}$ but replacing $n$ with $n-2$ and $k$ with $k-1$ gives $$f^{2^{k-2}1^{n-2k}} = \frac{(n-2)-2(k-1)+1}{k-1}\binom{n-2}{k-2} = \frac{n-2k+1}{n-1}\binom{n-1}{k-1}.$$ Putting this together, the multiplicity of $\mc{S}^\lam$ in $H_{1,0}(G)$ is $$(n-1)\cdot \frac{n-2k+1}{n-1}\binom{n-1}{k-1} - \frac{n-2k+1}k \binom{n}{k-1} = (n-2k+1)\left(\binom{n-1}{k-1}-\frac1k\binom{n}{k-1}\right),$$ as claimed. \qed

\section{Proof of the conjectured values from Table \ref{tab:homology}}\label{sec:table_vals}

We reproduce Table \ref{tab:homology} here for convenience:
\begin{center}
    \renewcommand{\arraystretch}{1.5}
    \begin{tabular}{|c|c|}
        \hline
        $\B{G}$ & $\B{H_{1,0}(G;\mb{C})}$ \\ \hline
        \im{figs/star4} & $\mc{S}_{2^2}$ \\ \hline
        \im{figs/star5} & $\mc{S}_{2^21}^{\oplus3} \oplus \mc{S}_{32}$ \\ \hline
        \im{figs/star6} & $\mc{S}_{2^21^2}^{\oplus6} \oplus \mc{S}_{2^3}^{\oplus5} \oplus \mc{S}_{321}^{\oplus4} \oplus \mc{S}_{42}$ \\ \hline
        \im{figs/star7} & $\mc{S}_{2^21^3}^{\oplus10} \oplus \mc{S}_{2^31}^{\oplus16} \oplus \mc{S}_{321^2}^{\oplus10} \oplus \mc{S}_{32^2}^{\oplus 9} \oplus \mc{S}_{421}^{\oplus5} \oplus \mc{S}_{52}$ \\ \hline
    \end{tabular}
\end{center}
Corollary \ref{cor:l21^(n-l-2)} immediately implies the following multiplicities from Table \ref{tab:homology}:
\begin{center}
    \renewcommand{\arraystretch}{2.2}
    \begin{tabular}{|c|c|}
        \hline
        \tb{Specht module} & \tb{Multiplicity} \\ \hline
        $\mc{S}_{2^2}$ in $H_{1,0}\left(\im{figs/star4}\right)$ & $\dbinom{4-2}2 = \dbinom22 = 1$ \\ \hline
        $\mc{S}_{2^21}$ in $H_{1,0}\left(\im{figs/star5}\right)$ & $\dbinom{5-2}2 = \dbinom32 = 3$ \\ \hline
        $\mc{S}_{32}$ in $H_{1,0}\left(\im{figs/star5}\right)$ & $\dbinom{5-2}3 = \dbinom33 = 1$ \\ \hline
        $\mc{S}_{2^21^2}$ in $H_{1,0}\left(\im{figs/star6}\right)$ & $\dbinom{6-2}2 = \dbinom42 = 6$ \\ \hline
        $\mc{S}_{321}$ in $H_{1,0}\left(\im{figs/star6}\right)$ & $\dbinom{6-2}3 = \dbinom43 = 4$ \\ \hline
        $\mc{S}_{42}$ in $H_{1,0}\left(\im{figs/star6}\right)$ & $\dbinom{6-2}4 = \dbinom44 = 1$ \\ \hline
        $\mc{S}_{2^21^3}$ in $H_{1,0}\left(\im{figs/star7}\right)$ & $\dbinom{7-2}2 = \dbinom52 = 10$ \\ \hline
        $\mc{S}_{321^2}$ in $H_{1,0}\left(\im{figs/star7}\right)$ & $\dbinom{7-2}3 = \dbinom53 = 10$ \\ \hline
        $\mc{S}_{421}$ in $H_{1,0}\left(\im{figs/star7}\right)$ & $\dbinom{7-2}4 = \dbinom54 = 5$ \\ \hline
        $\mc{S}_{52}$ in $H_{1,0}\left(\im{figs/star7}\right)$ & $\dbinom{7-2}5 = \dbinom55 = 1$ \\ \hline
    \end{tabular}
\end{center}
Corollary \ref{cor:2^k1^(n-2k)} then implies two of the others:
\begin{center}
    \renewcommand{\arraystretch}{2.2}
    \begin{tabular}{|c|c|}
        \hline
        \tb{Specht module} & \tb{Multiplicity} \\ \hline
        $\mc{S}_{2^3}$ in $H_{1,0}\left(\im{figs/star6}\right)$ & $(6-2\cdot3+1)\left(\dbinom{6-1}{3-1} - \dfrac13\dbinom6{3-1}\right) = \dbinom52 - \dfrac13\dbinom62 = 5$ \\ \hline
        $\mc{S}_{2^31}$ in $H_{1,0}\left(\im{figs/star7}\right)$ & $(7-2\cdot3+1)\left(\dbinom{7-1}{3-1} - \dfrac13\dbinom7{3-1}\right) = 2\left(\dbinom62 - \dfrac13\dbinom72\right) = 16$ \\ \hline
    \end{tabular}
\end{center}
The only other multiplicity to check is $\mc{S}_{32^2}^9$ in $H_{1,0}\left(\im{figs/star7}\right)$. For that one, we directly use the formula $$\binom{n-1}{\ell-1}f^{2^k1^{n-\ell-2k}} - f^\lam$$ from Theorem \ref{thm:general_l_k}, which in this case becomes $$\binom{7-1}{3-1}f^{2^2} - f^{32^2} = \binom62 f^{2^2} - f^{32^2} = 15f^{2^2} - f^{32^2}.$$ We have $f^{2^2} = 2,$ since the two SYT of shape $2^2$ are $$\young(12,34)\hspace{1cm}\tn{and}\hspace{1cm}\young(13,24).$$ For $f^{32^2},$ we use the hook length formula. Below, we label each of the 7 boxes with the length of its hook:
$$\young(541,32,21).$$ Thus, the hook length formula gives $$f^{32^2} = \frac{7!}{5\cdot4\cdot3\cdot2\cdot2\cdot1\cdot1} = 21.$$ Plugging these values in gives a multiplicity of $$15\cdot2 - 21 = 9,$$ as claimed in the table.

For the stars on 4 through 7 vertices, it follows by dimension counting that the full decompositions of $H_{1,0}(G)$ are as claimed in Table \ref{tab:homology}, because the authors of \cite{chandler2023strength} computed what the dimension of $H_{1,0}(G)$ is in each case. \qed

\section*{Acknowledgements}

The author thanks Oliver Pechenik for telling her about the chromatic symmetric homology, and she thanks Oliver Pechenik, Karen Yeats, and Sophie Spirkl for providing helpful comments. She was partially supported by the Natural Sciences and Engineering Research Council of Canada (NSERC) grant RGPIN-2022-03093.

\printbibliography

\end{document}